\documentclass[journal]{IEEEtran}

\usepackage{etex}
\usepackage[vlined,ruled]{algorithm2e}
\usepackage[dvipsnames]{xcolor}
\usepackage{pgfpages}
\usepackage[cmex10]{amsmath}		
\usepackage{amssymb, mathrsfs}
\usepackage{cite}

\usepackage[plainpages=false]{hyperref}
\hypersetup{
    colorlinks=true,       
    linkcolor=black,        
    citecolor=black,        
    breaklinks=true,
    urlcolor=blue,
}

\usepackage{dsfont}
\usepackage{siunitx}
\usepackage{verbatim}									
\usepackage{mathtools}									
\usepackage{siunitx}									
\usepackage{mdframed}

\usepackage{graphicx}
\usepackage{epstopdf}

\usepackage[font=small]{caption}
\usepackage{subcaption}
\usepackage[activate={true,nocompatibility},final,tracking=true,kerning=true,spacing=true,factor=1100,stretch=10,shrink=10]{microtype}

\usepackage{array}
\usepackage{multirow}
\usepackage{enumerate}
\usepackage{booktabs}

\graphicspath{{./fig/}}

\usepackage{pgfplots}
\pgfplotsset{width=8cm,compat=newest}
\newcommand*{\damping}{0.006}%
\newcommand*{\freq}{25}%
\pgfmathsetmacro{\freqd}{sqrt(1-(\damping)^2)*\freq}%

\pgfplotsset{
    standard/.style={
    axis x line=middle,
    axis y line=middle,
    enlarge x limits=0.15,
	enlarge y limits=0.15,
	every axis plot post/.style={mark options={fill=black}},
	}
}

\pgfplotsset{%
    ,compat=1.12
    ,every axis x label/.style={at={(current axis.right of origin)},anchor=north west}
    ,every axis y label/.style={at={(current axis.above origin)},anchor=north east}
    }

\usepackage{tikz}
\usetikzlibrary{arrows}
\usetikzlibrary{decorations.pathmorphing}
\usetikzlibrary{decorations.markings}
\usetikzlibrary{snakes}
\usetikzlibrary{matrix}
\usetikzlibrary{calc}
\usetikzlibrary{patterns}
\usetikzlibrary{positioning}
\usetikzlibrary{shapes}
\usetikzlibrary{shapes.symbols}
\usetikzlibrary{shapes.geometric}
\usetikzlibrary{shadows.blur}
\usetikzlibrary{shapes.arrows}
\usetikzlibrary{shapes.multipart}
\usetikzlibrary{shapes.callouts}
\usetikzlibrary{shapes.misc}

\tikzstyle{every node}=[font=\small]
\tikzstyle{every path}=[line width=0.8pt,line cap=round,line join=round]

\usepackage[american,smartlabels]{circuitikz}
\ctikzset{bipoles/thickness=1}
\ctikzset{bipoles/length=0.8cm}
\ctikzset{bipoles/diode/height=.375}
\ctikzset{bipoles/diode/width=.3}
\ctikzset{tripoles/thyristor/height=.8}
\ctikzset{tripoles/thyristor/width=1}
\ctikzset{bipoles/vsourceam/height/.initial=.7}
\ctikzset{bipoles/vsourceam/width/.initial=.7}

\newcommand{\real}{\mathbb{R}}

\newcommand{\setdef}[2]{\{#1 \;|\; #2\}}

\newcommand{\vect}[1]{\mathbbold{#1}}
\newcommand{\vones}[1][]{\vect{1}_{#1}}

\DeclareSymbolFont{bbold}{U}{bbold}{m}{n}
\DeclareSymbolFontAlphabet{\mathbbold}{bbold}

\newcommand{\map}[3]{#1: #2 \rightarrow #3}

\newcommand{\T}{\mathsf{T}} 

\newcommand\oprocendsymbol{\hbox{$\square$}}
\newcommand\oprocend{\relax\ifmmode\else\unskip\hfill\fi\oprocendsymbol}

\newtheorem{theorem}{Theorem}[section]
\newtheorem{lemma}[theorem]{Lemma}

\newtheorem{assumption}[theorem]{Assumption}

\newenvironment{pfof}[1]{\vspace{1ex}\noindent{\itshape Proof of
    #1:}\hspace{0.5em}} {\hfill\oprocend\vspace{1ex}}

\makeatletter
\newcommand{\vast}{\bBigg@{4}}
\newcommand{\Vast}{\bBigg@{5}}

\newcommand{\define}{\ensuremath{\triangleq}}

\usepackage{enumitem}

\renewenvironment{quote}{%
   \list{}{%
     \leftmargin0.5cm   
     \rightmargin\leftmargin
   }
   \item\relax
}
{\endlist}

\usepackage{enumerate}
\usepackage{cuted}
\usepackage{arydshln}
\newcommand{\NI}{\mathsf{NI}}

\begin{document}
%
\title{A Dynamic Stability and Performance Analysis of Automatic Generation Control}
%
%
%

\author{John~W.~Simpson-Porco,~\IEEEmembership{Member,~IEEE}
\thanks{J.~W.~Simpson-Porco is with the Department of Electrical and Computer Engineering, University of Toronto, 10 King's College Road,
Toronto, ON, M5S 3G4, Canada. Email: {\tt jwsimpson@ece.utoronto.ca}.}
}

%
%

\markboth{Submitted for publication. This version: \today}%
{Shell \MakeLowercase{\textit{et al.}}: Bare Demo of IEEEtran.cls for Journals}
%



\maketitle

\begin{abstract}
Automatic generation control (AGC) is one of the most important coordinated control systems present in modern interconnected power systems. Despite being heavily studied, no interconnected dynamic stability and performance analysis of AGC is available in the literature. This paper presents such an analysis for a class of multi-area interconnected nonlinear power systems, providing a nonlinear stability proof and examining dynamic performance from the perspectives of non-interactive tuning, response speed, and sensitivity to disturbances. The key insight is that dynamic stability and performance can be rigorously captured by a reduced dynamic model, which depends only on the AGC controller parameters and on the area frequency characteristic constants. Our analytical results clarify some of the historical controversy concerning the tuning of frequency bias constants in AGC for dynamic stability and performance, and are validated by simulations on a detailed test system. 
\end{abstract}



\section{Introduction}
\label{Sec:Introduction}



Modern large-scale AC power systems consist of interconnections between areas managed by local operators called balancing authorities (BAs). Mismatch between generation and load within each area is compensated through a hierarchy of control layers operating on different spatial and temporal scales. The ``secondary'' layer of control, called Automatic Generation Control (AGC), acts as an interface between the speed governing (``primary'') controllers of generation units and the global (``tertiary'') economic decisions computed via optimization. The purpose of AGC as an online controller is continuously reallocate generation to eliminate generation-load mismatch \emph{within} each area, subject to security and economy of operation. This estimation and allocation procedure ensures that the system frequency and all net inter-area power exchanges remain at their scheduled values. 



The distributed biased net-interchange concept, independently put forward by N. Cohn and R. Brandt \cite{NC:83}, allows each BA to achieve secondary control objectives using only local measurements by constructing an area control error (ACE) and regulating it to zero. Successfully deployed since the late 1940's, AGC has a long and extensive history of study, and its evolution continues to be a topic of considerable academic and industrial interest. We make no attempt to summarize the historical literature here; the clearest textbook treatment available can be found in \cite{AJW-BFW:96}, with a similar treatment in \cite{PK:94}, and we mention in particular \cite{FPD-RJM-WFB:73a,FPD-RJM-WFB:73b,HGK-KCK-AB:75,WBG:78,HG-JS:80,JC:85,MSC:86,NJ-LSV-DNE-LHF-AGH:92,IEEE-Report-AGC:79} as a selection of outstanding historical and/or practitioner perspectives on AGC.

Industry implementations of AGC are incredibly varied, and often include logical subroutines for reducing the activity of the AGC system and for handling system-specific conditions and operator preferences (e.g., \cite{FPD-RJM-WFB:73b,CWT-RLC:76}). The fundamental underlying characteristic that enables the success of these diverse implementations is the slow speed of operation of AGC (typically, minutes) relative to primary control dynamics (seconds to tens of seconds). This slow speed is \emph{intentional}, as the goal of AGC is smooth re-balancing of generation and load inbetween economic re-dispatch. The slow speed is also \emph{necessary}: models of primary frequency dynamics (including energy conversion, turbine-governor, and load dynamics) are subject to considerable uncertainty, and sampling/communication/filtering processes introduce unavoidable delays. To ensure closed-loop stability, all practical AGC systems must be sufficiently slow so that no significant dynamic interaction occurs between the integral loop and the primary frequency dynamics. As put succinctly in \cite{NJ-LSV-DNE-LHF-AGH:92}, ``\emph{reduction in the response time of AGC is neither possible nor desired}'' and ``\emph{attempting to do so serves no particular economic or control purposes.}'' 

Our work in this paper begins with the observation that the literature lacks a simple and compelling \emph{dynamic} stability and performance analysis of even elementary AGC implementations. The standard textbook analyses \cite{AJW-BFW:96,PK:94} are based on equilibrium analysis, and do not consider dynamic convergence to the equilibrium. Notable exceptions are the treatments in \cite{MDI-SL:96, MI-JZ:00}, but these analyses focus on reduced-order area models, and do not analyze the interconnected dynamics of AGC systems. Academic papers have tended to assess dynamic stability/performance of AGC via simulation, using detailed dynamic models of the area inertial and primary responses. Unfortunately, the results of such simulation studies are always highly dependent on the underlying modelling assumptions or test system, as each power system is unique in its dynamic response to disturbances. This has resulted in a great deal of confusion and controversy in the literature regarding ``optimal'' tuning of frequency bias factors $b$, dating to the (infamous) paper \cite{OIE-CEF:70}. NERC/ENTSO-E standards and nearly all literature specify that $b$ should be tuned equal to the frequency response characteristic (FRC) $\beta$ of the control area, with a preference towards overbiasing.\footnote{The determination of the FRC constant $\beta$ is itself a fraught exercise. Indeed, there is no one constant value of $\beta$ for any real power system, as it varies seasonally with load composition, varies with unit commitment and dispatch point, and varies dynamically depending on which governors are operating within their frequency deadbands. These issues are well outside our scope; see \cite{NJ-LSV-DNE-LHF-AGH:92,NERC:11,AR:11,MS-EI-AR-GA:12} for further reading.} This tuning has its origins in a static analysis due to Cohn \cite{NC:56}, but has no clear basis in terms of dynamic performance. One of our goals in this paper is to scrutinize this wisdom with a dynamic analysis.

\emph{Current trends in advanced frequency control:} A related line of research \textemdash{} somewhat disjoint from the industry implementations of AGC \textemdash{} considers the application of advanced control techniques for secondary control. An early reference is \cite{CEF-OIE:70}, and several surveys \cite{PKI-DPK:05,HHA-MEHG-RZ-EHF-PS:18,FD-SB-JWSP-SG:19a} are available summarizing aspects of this line of research, but  implementations of such advanced control methods are apparently rare.\footnote{Curiously, it was suggested as early as 1978 \cite{ED-NP:78} that advanced decentralized control was unlikely to offer significant advantages over the traditional AGC controller; this conclusion appears to have stood the test of time.} There is currently interest in understanding the limitations of AGC in the presence of renewable energy integration, and in subsequently developing modernized frequency control methods. Work in this direction includes model-predictive AGC \cite{ANV-IAH-JBR-SJW:08,AM-LI-KU:16,PRBM-PT:17}, various ``enhanced'' versions of AGC \cite{QL-MDI:12,DA:14,CL-JHC:17}, online gradient-type methods \cite{NL-CZ-LC:16, FD-SG:17, CZ-EM-SHL-JB:18}, and frequency or dynamics-aware dispatch and AGC \cite{DA-PWS-ADD:16,AAT-FZ-LX:11,GZ-JM-QW:19,PC-SD-YCC-MP:20}. While our work here is not directly motivated by these modernization topics, our technical approach will prove useful to other researchers in this area, and has already been applied for analysis of some advanced frequency control methods in \cite{JWSP:20a,JWSP:20b}.




\smallskip

\emph{Contributions:} This paper presents a formal dynamic stability and performance analysis of AGC, considering two conventional implementations which incorporate limiters for the control signals. The distinguishing characteristic of our analysis is that it accepts the necessity of time-scale separation between AGC and primary control, and rigorously follows this to its logical analysis conclusion using singular perturbation theory \cite{HKK:02}. As a result, the details of the inertial/primary control dynamics are rendered less important, which allows us to assess the properties AGC which are attributable purely to the decentralized control structure and the essential physics of power systems. We prove that when implemented on an interconnected nonlinear power system, both conventional AGC implementations are stable for any tuning of bias factors and time constants. We assess performance of AGC from several perspectives, including tuning for area-wise non-interaction, eigenvalue analysis, and sensitivity to disturbances. Our analysis makes explicit how the dynamic performance of AGC varies with the choices of bias factors and the FRC's of the areas, and \textemdash{} with some caveats \textemdash{} validates the orthodox tuning rule that overbiasing should be preferred to underbiasing.
 
Aside from fundamental theoretical interest, our dynamic analysis provides some practical insights into the tuning of bias factors for AGC. One interesting insight is that \textemdash{} for the textbook AGC implementation in \cite{AJW-BFW:96,PK:94} \textemdash{} the conventional wisdom of setting the AGC bias equal to the area FRC \emph{does not} lead to non-interacting behaviour between areas, and in fact, no such ideal bias setting exists for the textbook implementation.

\smallskip

\emph{Paper Organization:} Section \ref{Sec:Prelim} lays out our assumptions on the power system models and AGC implementations under consideration. Section \ref{Sec:AGCStability} contains our main nonlinear stability analysis results. Section \ref{Sec:AGCPerformance} studies performance of AGC as a function of bias tuning from several perspectives. Section \ref{Sec:Simulations} provides simulation results supporting our main theoretical results. The paper concludes with a brief discussion (Section \ref{Sec:Discussion}) and points to future research avenues (Section \ref{Sec:Conclusions}).

\smallskip

\emph{Notation:} The $N \times N$ identity matrix is $I_N$. Given a vector $x \in \real^N$, $x^{\T}$ denotes its transpose. 
Throughout, $\vones[N]$ is $N$-dimensional vector with unit entries; we will usually omit the subscript, and use the compact notation $\vones[]^{\T}x = \sum_{i=1}^{N}x_i$ for the sum of the elements of a vector $x \in \real^N$. Given a collection of numbers or column vectors $(v_1,\ldots,v_N)$, we let $\mathrm{col}(v_1,\ldots,v_N)$ denote the concatenated column vector and $\mathrm{diag}(v_1,\ldots,v_N)$ denote the associated block diagonal matrix. For a square matrix $A$ with real eigenvalues $\lambda_i(A)$, we order them as $\lambda_{1}(A) \leq \cdots \leq \lambda_{N}(A)$.


\section{Nonlinear Power System Model and Automatic Generation Control}
\label{Sec:Prelim}


\subsection{Interconnected Power System Model}
\label{Sec:PowerModel}

We consider an interconnected power system consisting of $N$ areas, and label the set of areas as $\mathcal{A} = \{1,\ldots,N\}$. In area $k \in \mathcal{A}$, suppose there are $m_k$ generators; we label the set of generators as $\mathcal{G}_{k} = \{1,\ldots,m_k\}$, and we let $\mathcal{G}_{k}^{\rm AGC} \subseteq \mathcal{G}_k$ denote the subset of generators which participate in AGC. Each generator $i \in \mathcal{G}_{k}$ has an electrical power output $P_{k,i}$ and its governor accepts a load reference command $u_{k,i}$, with $u_{k,i}^{\star}$ denoting the base reference determined by economic dispatch. Load reference commands are restricted to the limits $\underline{u}_{k,i} \leq u_{k,i} \leq \overline{u}_{k,i}$. If generator $i \in \mathcal{G}_k$ does not participate in AGC, then $u_{k,i}$ is fixed at $u_{k,i}^{\star} \in [\underline{u}_{k,i},\overline{u}_{k,i}]$. For area $k$ we define the vector variables $P_k = \mathrm{col}(P_{k,1},\ldots,P_{k,m_k})$, with $u_k$ and $u_k^{\star}$ defined similarly. We let $\Delta f_k = f_k - f_k^{\star}$ and $\Delta \NI_{k} = \NI_{k} - \NI_{k}^{\star}$ be measurements of frequency and net interchange deviation for area $k$, with respect to their scheduled values.\footnote{Note that $\Delta f_k$ could be a weighted average of measurements internal to the area; this will have no effect on our analysis. For the net interchange deviation, a positive value corresponds to a net flow out of the area.}

We collect all variables for \emph{all} areas into larger stacked vectors $P = \mathrm{col}(P_1,\ldots,P_N)$, with $u \in \real^{m_1+\cdots+m_N}$, $\Delta f \in \real^N$ and $\Delta \NI \in \real^N$ similarly defined. The entire interconnected power system is described by a nonlinear ODE model
\begin{subequations}\label{Eq:NonlinearPowerSystem}
\begin{align}
\label{Eq:NonlinearPowerSystem-1}
\dot{x}(t) &= F(x(t),u(t),w(t))\\
\label{Eq:NonlinearPowerSystem-2}
\mathrm{col}(\Delta f(t),\Delta \NI(t),P(t)) &= h(x(t),u(t),w(t)),
\end{align}
\end{subequations}
where $x(t)$ is the vector of states and $F,h$ are appropriate functions. The dynamics \eqref{Eq:NonlinearPowerSystem-1} are assumed to already include the typical filters for the measurements specified in \eqref{Eq:NonlinearPowerSystem-2}, and and \eqref{Eq:NonlinearPowerSystem} may have been obtained from a more general differential-algebraic model under appropriate regularity conditions \cite{DJH-IMYM-90}. The disturbance $w(t)$ can model changes to the frequency/net interchange schedules and reference changes to other control loops (e.g., AVRs), but most importantly $w$ includes the unmeasured net load deviation $\Delta P_{k}^{\rm L}$ for each area $k \in \mathcal{A}$. The precise model \eqref{Eq:NonlinearPowerSystem} is never known in reality, and we make no attempt to specify it. We will instead assume that the model satisfies some basic stability and steady-state properties.


\smallskip

\begin{assumption}[\bf Stable Nonlinear Power System Model]\label{Ass:PowerSystem}
There exist domains $\mathcal{X}, \mathcal{U},\mathcal{W}$ for $x$, $u$, and $w$ such that
\begin{enumerate}[label=(A\arabic*)]
\item \label{Ass:PowerSystem-1} \textbf{Model Regularity:} $F$, $h$, and their Jacobians are Lipschitz continuous on $\mathcal{X}$ uniformly in $(u,w) \in \mathcal{U} \times \mathcal{W}$;
\item \label{Ass:PowerSystem-2} \textbf{Existence of Steady-State:} there exists a continuously differentiable equilibrium map $\map{x_{\rm ss}}{\mathcal{U}\times\mathcal{W}}{\mathcal{X}}$ which is Lipschitz continuous on $\mathcal{U}\times\mathcal{W}$ and satisfies $0 = F(x_{\rm ss}(u,w),u,w)$ for all $(u,w) \in \mathcal{U} \times \mathcal{W}$;
\item \label{Ass:PowerSystem-3} \textbf{Uniform Exponential Stability of the Steady-State:} the steady-state $x_{\rm ss}(u,w)$ is locally exponentially stable, uniformly in the inputs $(u,w) \in \mathcal{U} \times \mathcal{W}$;
%
%
\item \label{Ass:PowerSystem-4} \textbf{Steady-State Synchronism and Interchange Balance:} for each $(u,w) \in \mathcal{U}\times\mathcal{W}$ the steady-state values of $(\Delta f,\Delta \NI ,P)$ determined by $(\Delta f,\Delta \NI,P) = h(x_{\rm ss}(u,w),u,w)$ satisfy the synchronism condition
\begin{equation}\label{Eq:Synchronous}
\Delta f_1 = \Delta f_2 = \cdots = \Delta f_N = \Delta f_{\rm ss}
\end{equation}
and the net interchange balance condition
\begin{equation}\label{Eq:NetBalance}
0 = \sum_{k\in\mathcal{A}}\nolimits \Delta \NI_k.
\end{equation}
\item \label{Ass:PowerSystem-5} \textbf{Area Balance and Lossless Tie Lines:}
the steady-state values $(\Delta f,\Delta \NI,P)$ from \ref{Ass:PowerSystem-4} additionally satisfy the area-wise balance conditions
\begin{subequations}\label{Eq:SteadyStateFormulas}
\begin{align}
\label{Eq:Balancek}
\sum_{i\in\mathcal{G}_k}\nolimits (P_{k,i}-u_{k,i}^{\star}) &= D_{k}\Delta f_{k} + \Delta P_{k}^{\rm L} + \Delta p_{k}^{\rm out}\\
\label{Eq:GovPower}
P_{k,i} &=  u_{k,i} - \tfrac{1}{R_{k,i}}\Delta f_k
\end{align}
\end{subequations}
for each $k \in \mathcal{A}$ and $i \in \mathcal{G}_{k}$,  where $D_{k} > 0$ models aggregate area load damping and $R_{k,i} > 0$ is the primary control gain of generator $i \in \mathcal{G}_{k}$. Moreover, all inter-area tie lines are lossless, which implies that the change in net tie-line power flow $\Delta p_{k}^{\rm out}$ out of each area $k$ satisfies
\begin{equation}\label{Eq:PowerFlowNI}
\Delta p_{k}^{\rm out} = \Delta\NI_{k}.
\end{equation}
\end{enumerate}
\end{assumption}

\medskip

Assumptions \ref{Ass:PowerSystem-1}--\ref{Ass:PowerSystem-3} say that the model is sufficiently regular, and that for any constant inputs $(u,w) \in \mathcal{U}\times\mathcal{W}$, there is a unique exponentially stable equilibrium state $x = x_{\rm ss}(u,w)$ which lives in some subset $\mathcal{X}$ of the normal operating region. In \ref{Ass:PowerSystem-4}, the synchronization condition \eqref{Eq:Synchronous} will always hold under normal operating conditions, and the interchange balance condition \eqref{Eq:NetBalance} always holds, as each tie-line is metered at a common point \cite{NERC:10bal}. The key model assumptions are in \ref{Ass:PowerSystem-5}. The area-wise balance conditions \eqref{Eq:SteadyStateFormulas} model area power balance and linear primary control, with \eqref{Eq:PowerFlowNI} being implied by lossless tie lines. Consideration of turbine-governor deadbands and tie line losses is outside our scope; see also Section \ref{Sec:Discussion}. To avoid analysis of saturated operation and integrator anti-windup implementations, we assume there is sufficient regulation capacity in each area.

\smallskip

\begin{assumption}[\bf Strict Local Feasibility]\label{Ass:Feasibility}
Each area has sufficient regulation capacity to meet the disturbance, i.e.,
\[
\Delta P^{\rm L}_{k} \in \mathcal{C}_k \define \Bigg(\sum_{i\in\mathcal{G}_k^{\rm AGC}} (\underline{u}_{k,i}-u_{k,i}^{\star}), \sum_{i\in\mathcal{G}_k^{\rm AGC}}(\overline{u}_{k,i}-u_{k,i}^{\star})\Bigg)
\]
for each area $k \in \mathcal{A}$, and we let $\mathcal{C} = \mathcal{C}_1 \times \cdots \times \mathcal{C}_N$.
\end{assumption}

\smallskip

By combining the formulas \eqref{Eq:Synchronous}--\eqref{Eq:PowerFlowNI}, one can calculate that the steady-state net interchange deviation $\Delta \NI_{k}$ for area $k$ and the steady-state frequency deviation $\Delta f_{\rm ss}$ are given by
\begin{subequations}
\begin{align} \label{Eq:SteadyStatePowerFlow}
&\begin{aligned} 
\hspace{-3em}\Delta \NI_{k} &= \tfrac{\beta-\beta_k}{\beta}(\vones[]^{\T}(u_k-u_k^{\star}) - \Delta P_{k}^{\rm L})\\
&\qquad - \tfrac{\beta_k}{\beta}\sum_{j\in\mathcal{A}\setminus\{k\}}\nolimits (\vones[]^{\T}(u_j-u_j^{\star}) - \Delta P_{j}^{\rm L})
\end{aligned}\\
\label{Eq:SteadyStateFrequency}
\Delta f_{\rm ss} &= \tfrac{1}{\beta}\sum_{k\in\mathcal{A}}\nolimits (\vones[]^{\T}(u_k-u_k^{\star}) - \Delta P_{k}^{\rm L}),
\end{align}
\end{subequations}
%
%
%
where $\beta_k = D_k + R^{-1}_k$ is the \emph{frequency response characteristic} (FRC) of area $k \in \mathcal{A}$ and $\beta = \sum_{k\in\mathcal{A}}\beta_k$ is the FRC of the interconnected system. For future use, note the identity
\begin{equation}\label{Eq:GovIdentity}
\vones[]^{\T}(u_k-u_k^{\star}) = \sum_{j\in\mathcal{G}_{k}^{\rm AGC}}\nolimits (u_{k,j}-u_{k,j}^{\star}).
\end{equation}

\subsection{Area Control Error}
\label{Sec:ACE}

Modulo signs and scalings, we use the standard NERC definition of the \emph{area control error} $\mathsf{ACE}_k$ for area $k \in \mathcal{A}$, which combines the measurements $\Delta \NI_{k}(t)$ and $\Delta f_{k}(t)$ as
\begin{equation}\label{Eq:ACE}
\begin{aligned}
\mathsf{ACE}_{k}(t) &\define \Delta \NI_{k}(t) + b_k\Delta f_k(t),\\
\end{aligned}
\end{equation}
where $b_k > 0$ is the \emph{frequency bias} for area $k$. The following result is straightforward to prove (Appendix \ref{App:1}) by combining \eqref{Eq:Synchronous}--\eqref{Eq:ACE}.

\smallskip

\begin{lemma}[\bf Steady-State Zeroing of ACE]\label{Lem:ACE}
If the interconnected power system \eqref{Eq:NonlinearPowerSystem} is in steady-state as specified in Assumption \ref{Ass:PowerSystem}, then the following statements are equivalent: 
\begin{enumerate}[label=(\roman*)]
\item \label{Lem:ACE-1} $\sum_{i\in\mathcal{G}_{k}^{\rm AGC}}\nolimits (u_{k,i}-u_{k,i}^{\star}) = \Delta P_{k}^{\rm L}$ for all areas $k \in \mathcal{A}$;
\item \label{Lem:ACE-2} $\Delta f_{\rm ss} = 0$ and $\Delta \NI_{k} = 0$ for all areas $k \in \mathcal{A}$;
\item \label{Lem:ACE-3} $\mathsf{ACE}_k = 0$ for all areas $k \in \mathcal{A}$.
\end{enumerate}
\end{lemma}

In other words, zeroing all ACEs is equivalent to frequency and net interchange regulation, which is in turn equivalent to local balancing of all net loads.


\subsection{Textbook and Simplified AGC Dynamics}
\label{Sec:AGC}

The simplest implementation of AGC one encounters in the literature integrates the ACE to produce an AGC control signal $\eta_k$ for area $k$ as
\begin{equation}\label{Eq:AGCSimple}
\tau_k \dot{\eta}_k(t) = -\mathsf{ACE}_k(t), \qquad k \in \mathcal{A},
\end{equation}
where $\tau_k > 0$ is the integral time constant, quoted in the literature as ranging from 30s up to 200s. This differs from the implementation found in standard textbooks, which includes an auxiliary feedback term involving electric power outputs of all generators within the area \cite{NC:60, AJW-BFW:96, DA-YCC-JZ-ADD-PWS:13, SVD-YCC-AAD-ADG:20}, formulated as
\begin{equation}\label{Eq:AGC}
\tau_{k} \dot{\eta}_{k}(t) = -\mathsf{ACE}_k(t) - \sum_{j\in\mathcal{G}_{k}}(u_{k,j}(t)-P_{k,j}(t)).
\end{equation}
Control actions from the AGC system are allocated across all participating generators $\mathcal{G}_{k}^{\rm AGC}$ such that their incremental costs of production (or, with lossess, delivery) are roughly equalized. A typical allocation rule including saturation is
\begin{equation}\label{Eq:Allocation}
u_{k,i} = \mathrm{sat}_{k,i}(u_{k,i}^{\star} + \alpha_{k,i}\eta_k), \quad i \in \mathcal{G}_k^{\rm AGC},
\end{equation}
where $\mathrm{sat}_{k,i}(v)$ saturates its argument to the limits $[\underline{u}_{k,i},\overline{u}_{k,i}]$, and $\{\alpha_{k,i}\}_{i\in\mathcal{G}_{k}^{\rm AGC}}$ are nonnegative \emph{participation factors} \cite{AJW-BFW:96} with normalization $\sum_{i\in\mathcal{G}_{k}^{\rm AGC}}\nolimits \alpha_{k,i} = 1$ for each area $k \in \mathcal{A}$.

\smallskip


\section{Dynamic Stability of Automatic Generation Control}
\label{Sec:AGCStability}

The closed-loop system consists of the interconnected power system \eqref{Eq:NonlinearPowerSystem} with either the textbook AGC controllers \eqref{Eq:AGC} or the simplified AGC controllers \eqref{Eq:AGCSimple}, and the allocation rules \eqref{Eq:Allocation}. Our analysis approach is based on a (rigorous) quasi steady-state analysis of \eqref{Eq:NonlinearPowerSystem}, where we assume the AGC dynamics are slow compared to the power system dynamics.  This analysis approach is strongly justified by the time-scale properties of AGC, as outlined in Section \ref{Sec:Introduction}. The key technical tool we employ is singular perturbation theory \cite{AS-HKK:84,HKK:02}, which rigorously justifies this approximation.

We begin examining the value the ACE takes when the power system is in steady-state (Assumption \ref{Ass:PowerSystem}). Substituting \eqref{Eq:SteadyStatePowerFlow} and \eqref{Eq:SteadyStateFrequency} into \eqref{Eq:ACE}, we obtain
\[
\begin{aligned}
\mathsf{ACE}_{k} &= \tfrac{\beta+b_k-\beta_k}{\beta}(\vones[]^{\T}(u_k-u_k^{\star}) - \Delta P_{k}^{\rm L})\\
&\qquad  + \tfrac{b_k-\beta_k}{\beta}\sum_{j \in \mathcal{A}\setminus\{k\}}\nolimits (\vones[]^{\T}(u_j-u_{j}^{\star}) - \Delta P_{j}^{\rm L}).
\end{aligned}
\]
Similarly, in steady-state we have from \eqref{Eq:GovPower},\eqref{Eq:SteadyStateFrequency} that
\[
\begin{aligned}
\vones[]^{\T}(u_{k}-P_{k}) &= \tfrac{1}{R_k}\Delta f_{\rm ss} = \tfrac{1}{R_k}\tfrac{1}{\beta}\sum_{j\in\mathcal{A}} (\vones[]^{\T}(u_j-u_j^{\star}) - \Delta P_{j}^{\rm L}).
\end{aligned}
\]
Using \eqref{Eq:Allocation}, we may compactly write
\[
\vones[]^{\T}(u_k - u_k^{\star}) = \underbrace{\sum_{i\in\mathcal{G}_k}\nolimits(\mathrm{sat}_{k,i}(u_{k,i}^{\star} + \alpha_{k,i}\eta_k) - u_{k,i}^{\star})}_{\define \varphi_{k}(\eta_k)}.
\]
Note that in the absence of saturation we simply have $\varphi_k(\eta_k) = \eta_k$. Substituting the above formulas into the textbook AGC controller \eqref{Eq:AGC} and finally writing everything in vector notation, we obtain the \emph{reduced textbook AGC dynamics}
\begin{equation}\label{Eq:SlowTimeScaleAGC}
\begin{aligned}
\tau \dot{\eta}
 &= -\mathcal{B}_{\rm txt}(\varphi(\eta)-\Delta P^{\rm L})\\
\mathsf{ACE} &= \mathcal{B}_{\sf ACE}(\varphi(\eta)-\Delta P^{\rm L}),
\end{aligned}
\end{equation}
where all variables are now stacked vectors, $\tau = \mathrm{diag}(\tau_1,\ldots,\tau_N)$, and the interconnection matrices $\mathcal{B}_{\rm txt}, \mathcal{B}_{\sf ACE} \in \real^{N \times N}$ are defined in \eqref{Eq:BMatrix}, \eqref{Eq:BSimpleMatrix}. The model \eqref{Eq:SlowTimeScaleAGC} defines a nonlinear dynamic system with state vector $\eta \in \real^N$, input vector $\Delta P^{\rm L} \in \real^N$, output vector $\mathsf{ACE} \in \real^N$. Note that the matrices $\mathcal{B}_{\rm txt}$ and $\mathcal{B}_{\sf ACE}$ in \eqref{Eq:BMatrix}, \eqref{Eq:BSimpleMatrix} depend \emph{only} on the frequency bias constants $b_k$, FRCs $\beta_k$, and the primary control settings $R_k$. An identical set of arguments can be made using the simplified AGC dynamics \eqref{Eq:AGCSimple}, and in that case one instead obtains the \emph{reduced simplified AGC dynamics}
\begin{equation}\label{Eq:SlowTimeScaleAGC2}
\begin{aligned}
\tau \dot{\eta} &= -\mathcal{B}_{\sf ACE}(\varphi(\eta)-\Delta P^{\rm L})\\
 \mathsf{ACE} &= \mathcal{B}_{\sf ACE}(\varphi(\eta)-\Delta P^{\rm L}).
 \end{aligned}
\end{equation}
The reduced dynamic equations accurately model the system evolution \emph{after} the action of primary control \cite[Chap. 11]{HKK:02}; we will illustrate this via simulation in Section \ref{Sec:Simulations}. Our main nonlinear stability result covers both implementations.

\smallskip

\begin{theorem}[\bf Closed-Loop Stability with Automatic Generation Control]\label{Thm:AGCStable}
Consider the interconnected power system \eqref{Eq:NonlinearPowerSystem} under Assumptions \ref{Ass:PowerSystem}--\ref{Ass:Feasibility}, with either the simplified AGC controllers \eqref{Eq:AGCSimple} or the textbook AGC controllers \eqref{Eq:AGC}, and with the allocation rule \eqref{Eq:Allocation}. If the smallest AGC time constant is sufficiently large, then the closed-loop system possesses a unique exponentially stable equilibrium point $(\bar{x},\bar{\eta}) \in \mathcal{X} \times \real^{N}$ and $\mathsf{ACE}_k(t) \to 0$ as $t \to \infty$ for all areas $k \in \mathcal{A}$.
\end{theorem}

\smallskip

Theorem \ref{Thm:AGCStable} states that \textemdash{} with the usual time-scales of operation \textemdash{} closed-loop stability of both textbook and simplified AGC systems is guaranteed for \emph{any} tuning of bias factors and time constants. Note that this ``unconditional stability'' is consistent with engineering practice, in which balancing authorities independently tune their AGC controllers without coordination. In the proof (Appendix \ref{App:1}), we show that assessing closed-loop stability essential boils down to assessing stability of the reduced nonlinear dynamics \eqref{Eq:SlowTimeScaleAGC}/\eqref{Eq:SlowTimeScaleAGC2}, which describe the AGC system behaviour on the long time-scale.

\begin{figure*}
\begin{equation}\label{Eq:BMatrix}
\mathcal{B}_{\rm txt} \define \frac{1}{\beta}\begin{bmatrix}
\beta + b_1 + R_1^{-1} - \beta_1 & b_1 + R_1^{-1} - \beta_1 & \cdots & b_1 + R_1^{-1} - \beta_1\\
b_2 + R_2^{-1} - \beta_2 & \beta + b_2 + R_2^{-1} - \beta_2 & \cdots & \cdots \\
\vdots & \cdots & \ddots & b_{N-1} + R_{N-1}^{-1} - \beta_{N-1}\\
b_N + R_N^{-1} - \beta_N & \cdots & b_N + R_N^{-1} - \beta_N & \beta + b_N + R_N^{-1} - \beta_N 
\end{bmatrix}
\end{equation}
\end{figure*}

\begin{figure*}
\begin{equation}\label{Eq:BSimpleMatrix}
\mathcal{B}_{\sf ACE} \define \frac{1}{\beta}\begin{bmatrix}
\beta + b_1 - \beta_1 & b_1 - \beta_1 & \cdots & b_1 - \beta_1\\
b_2 - \beta_2 & \beta + b_2 - \beta_2 & \cdots & \cdots \\
\vdots & \cdots & \ddots & b_{N-1} - \beta_{N-1}\\
b_N - \beta_N & \cdots & b_N  - \beta_N & \beta + b_N - \beta_N 
\end{bmatrix}
\end{equation}
\end{figure*}

\section{Dynamic Performance of Automatic Generation Control}
\label{Sec:AGCPerformance}

In Section \ref{Sec:AGCStability} we established that under mild assumptions on the power system dynamics, standard AGC implementations are provably stable. We now explore the implications of our results for tuning and dynamic performance. The question of \emph{how} to quantify performance of AGC systems is a complex and historically controversial one. In practice, a good performance measure involves a mix of technical and non-technical considerations; the older technical references \cite{FPD-RJM-WFB:73a,FPD-RJM-WFB:73b,HGK-KCK-AB:75,WBG:78,HG-JS:80,JC:85,MSC:86,NJ-LSV-DNE-LHF-AGH:92,IEEE-Report-AGC:79} in the introduction provide substantial discussion on this point.

We will focus on quantifying the dynamic \emph{control} performance of AGC systems, via several metrics which follow naturally from our time-scale separation approach in Section \ref{Sec:AGCStability}. Our analysis will reveal some of the fundamental performance limitations of AGC which arise \emph{purely} due to the decentralized control structure and the selection of frequency bias factors. These performance limitations are intrinsic to the control architecture: they are always present, and are \emph{entirely unrelated} to other practical factors which may additionally degrade performance, such as abnormal system operating conditions, filtering and communication delays, measurement sampling periods, and ramping limitations of the participating generators. In other words, even an \emph{ideal} AGC implementation will still be subject to the performance limitations we note next.


\subsection{Revising Classical Bias Tuning for Non-Interaction}
\label{Sec:BiasTuning}

The ``optimal'' choice of the bias setting for use in AGC has been a topic of substantial  interest and controversy, dating back to the 1950's. In \cite{NC:56}, Cohn argued \textemdash{} based on static equilibrium analysis of the ACE (see, e.g., \cite{PK:94}) \textemdash{} that each area should set its bias $b_k$ equal to its FRC $\beta_k$. In doing so, each area will minimally respond to disturbances occurring in other areas. This reccomendation has been widely accepted, and is standard for both NERC and ENTSO-E \cite{NERC:11,UCTE-App1:04}. We now scrutinize this recommendation in the context of our reduced dynamics \eqref{Eq:SlowTimeScaleAGC} and \eqref{Eq:SlowTimeScaleAGC2}, which model the interconnected dynamics on the time-scale after the action of primary control. 

\smallskip

\emph{Simple AGC Dynamics \eqref{Eq:SlowTimeScaleAGC2}:} If the bias $b_k$ in area $k$ is tuned such that $b_k = \beta_k$, then all off-diagonal elements in the $k$th row of $\mathcal{B}_{\sf ACE}$ become zero, and the $k$th row of the dynamics \eqref{Eq:SlowTimeScaleAGC2} simplifies to the simple single-input single-output system
\begin{equation}\label{Eq:Decoupled}
\tau_k \dot{\eta}_k = -\varphi_{k}(\eta_k) + \Delta P^{\rm L}_{k},\,\,\,\,\, \mathsf{ACE}_k = \varphi_k(\eta_k) - \Delta P^{\rm L}_{k}.
\end{equation}
This shows that the control signal for area $k$ converges with what is essentially a first-order response, and is not influenced by any other areas. If \emph{all} areas select $b_k=\beta_k$, then all AGC systems are \emph{non-interacting}; this provides a dynamic systems justification for Cohn's conclusion.

\smallskip

\emph{Textbook AGC dynamics \eqref{Eq:SlowTimeScaleAGC}:} Inspection of \eqref{Eq:SlowTimeScaleAGC} and \eqref{Eq:BMatrix} immediately shows that the tuning  $b_k = \beta_k$ does not lead us to the same conclusions as above. Indeed, one can quickly see by examining \eqref{Eq:BMatrix}, \eqref{Eq:BSimpleMatrix} that there does not exist \emph{any} tuning of bias factors which yields the non-interacting equations \eqref{Eq:Decoupled}. We claim this points to a fundamental design deficiency in the textbook AGC model \eqref{Eq:AGC}. A full investigation of this deficiency is beyond our scope here; we refer the reader to our related letter \cite{JWSP:20e}\footnote{\url{https://www.control.utoronto.ca/~jwsimpson/papers/2020e.pdf}} for further insights. In the remainder of this paper, we therefore focus on analyzing the simple AGC model \eqref{Eq:AGCSimple} with reduced dynamics \eqref{Eq:SlowTimeScaleAGC2}.



\subsection{Dynamic Performance}
\label{Sec:DynamicPerf}

In terms of tuning, Cohn argued in \cite{NC:56} that overbiasing (i.e., $b_k > \beta_k$) should be preferred to underbiasing (i.e., $b_k < \beta_k$). Writing later in \cite{NC:83}, the following \emph{dynamic} claim is made:

\begin{quote}
\emph{
``\ldots settings lower than the combined governor-load governing characteristics resulted in undesirable withdrawl of assistance to areas in need. Such withdrawl is appreciably greater in relative magnitude than additional assistance that would be provided if settings were above the combined governing characteristic.''
}
\end{quote}
The first part of the above statement is true and uncontroversial; we will however scrutinize the second part of the statement in the context of our reduced dynamics \eqref{Eq:SlowTimeScaleAGC2}. To focus in specifically on the effect of bias tuning, in the remainder of this section we will (i) ignore the effects of saturation, so that $\varphi(\eta) = \eta$, and (ii) assume equal time constants $\tau_k = \tau_\ell = \tau^{\prime}$ for some $\tau^{\prime} > 0$ and all $k,\ell \in \mathcal{A}$.
Under these assumptions \eqref{Eq:SlowTimeScaleAGC2} becomes the LTI system
\begin{equation}\label{Eq:SlowTimeScaleAGC2-LTI}
\begin{aligned}
 \dot{\eta} &= -\tfrac{1}{\tau^{\prime}}\mathcal{B}_{\sf ACE}(\eta-\Delta P^{\rm L})
 \end{aligned}
\end{equation}
with convergence rate governed by the eigenvalues of the matrix $-\tfrac{1}{\tau^{\prime}}\mathcal{B}_{\sf ACE}$. In vector notation, we can write 
\[
-\tfrac{1}{\tau^{\prime}}\mathcal{B}_{\sf ACE} = -\tfrac{1}{\tau^{\prime}}(I_N - \tfrac{1}{\beta}(\boldsymbol{\beta}-\boldsymbol{b})\vones[]^{\T}),
\]
where $\boldsymbol{b} = \mathrm{col}(b_1,\ldots,b_N)$ and $\boldsymbol{\beta} = \mathrm{col}(\beta_1,\ldots,\beta_N)$. It follows then from Lemma \ref{Lem:Inter} that $-\tfrac{1}{\tau^{\prime}}\mathcal{B}_{\sf ACE}$ has $N-1$ eigenvalues at $-\tfrac{1}{\tau^{\prime}}$, with its $N$th eigenvalue given by
\[
\lambda_{N} = -\tfrac{1}{\tau^{\prime}}\left(\sum_{k\in\mathcal{A}} \nolimits b_k\right)/\left(\sum_{k\in\mathcal{A}} \nolimits \beta_k\right) < 0.
\]
For underbiased tunings $b_k \leq \beta_k$ with strict inequality for at least one area, we have $\lambda_{N} > -1/\tau^{\prime}$, and therefore $\lambda_{N}$ becomes the dominant slow eigenvalue. For overbiased tunings where $b_k \geq \beta_k$ with strict inequality in at least one area, we have $\lambda_{N} < -1/\tau^{\prime}$, and hence the $N-1$ eigenvalues at $-1/\tau^{\prime}$ are dominant. This shows that underbiasing will lead to a system-wide response that is \emph{slower} than what one would expect based only on the time constant $\tau^{\prime}$ (see Figure \ref{Fig:KundurOverUnder} later).

Unfortunately, the eigenvalues do not provide information about the transient effect of overbiasing, nor are they informative about the dynamic responses of the AGC systems to disturbances. To adequately capture this information, we will examine transfer functions associated with the dynamics \eqref{Eq:SlowTimeScaleAGC2-LTI}.

\smallskip

\begin{figure}[ht!]
\begin{center}
\includegraphics[width=0.95\columnwidth]{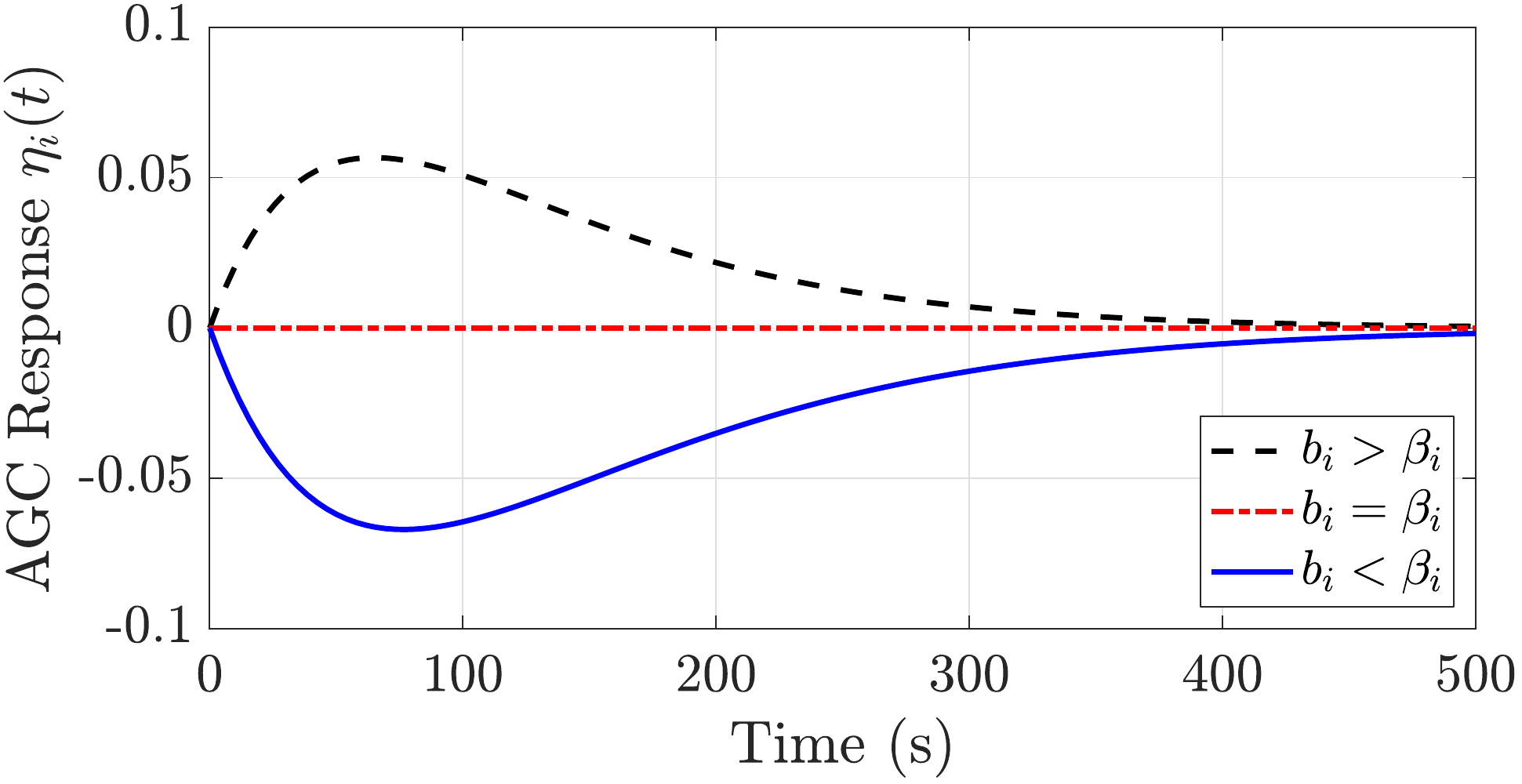}
\caption{Response of $\eta_i(t)$ due disturbance in area $j$; $\tau^{\prime} = 70$s.}
\label{Fig:OverUnderStep}
\end{center}
\end{figure}

\emph{Transfer function $T_{ij}(s)$ from net load disturbance $\Delta P_{j}^{\rm L}$ in area $j$ to AGC control signal $\eta_i$ in area $i$.} From \eqref{Eq:SlowTimeScaleAGC2-LTI}, the required transfer function is given by
\[
T_{ij}(s) = \frac{\eta_i(s)}{\Delta P_{j}^{\rm L}(s)} = \mathsf{e}_i^{\T}(\tau^{\prime} sI_{N}+\mathcal{B}_{\sf ACE})^{-1}\mathcal{B}_{\sf ACE}\mathsf{e}_j,
\]
where $\mathsf{e}_i \in \real^N$ denotes the $i$th unit vector in $\real^N$. Straightforward but tedious calculation using Lemma \ref{Lem:Inter} shows that
\[
\begin{aligned}
T_{ij}(s) 
%
%
%
&= \frac{1}{\tau^{\prime}s+1}\left[\delta_{ij} - \frac{1}{\beta}(\beta_i-b_i)\frac{\tau^{\prime}s}{\tau^{\prime}s+\tfrac{1}{\beta}\sum_{k}b_k}\right],
\end{aligned}
\]
where $\delta_{ij} = 1$ if $i = j$ and $0$ otherwise. We can use this transfer function to evaluate Cohn's claim. Consider applying a step load disturbance in area $j$, and examining the AGC response in area $i$; the generic response is plotted in Figure \ref{Fig:OverUnderStep} for (equally) underbiased and overbiased tunings of area $i$. 

\begin{figure}[ht!]
\begin{center}
\includegraphics[width=0.95\columnwidth]{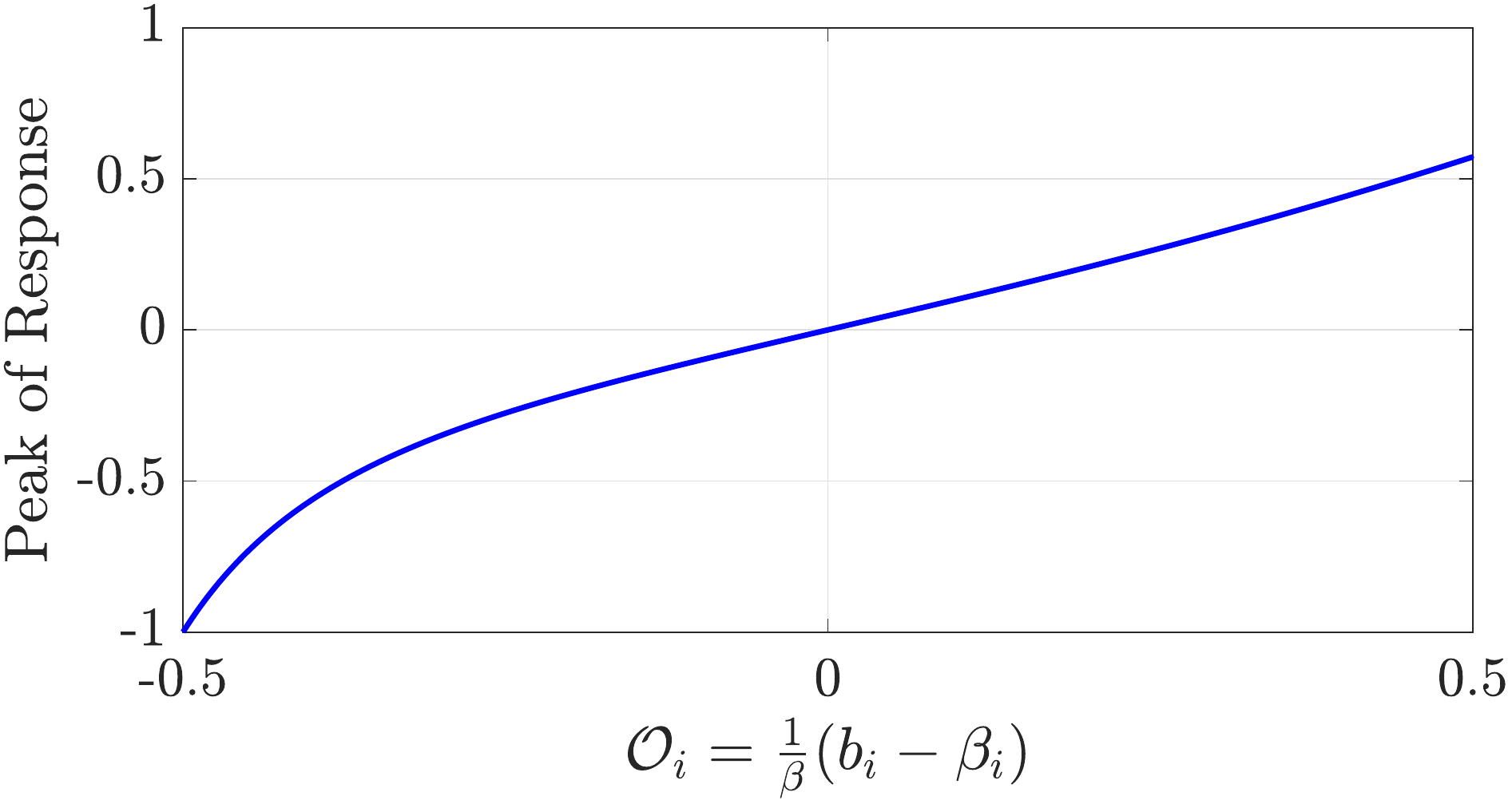}
\caption{Peak response of $\eta_i(t)$ due to under/overbiasing.}
\label{Fig:OverUnder}
\end{center}
\end{figure}

Note that the underbiased and overbiased responses are symmetric, but not perfectly so. To analytically quantify this, let $\mathcal{O}_i = (b_i-\beta_i)/\beta$ quantify the over/under biasing of area $i$ (note that due to division by the overall FRC $\beta$, $\mathcal{O}_i$ will be quite small in reality). Straightforward but cumbersome calculations show that the \emph{peak value} of the response $\eta_i(t)$ in area $i$ is given by
\[
\mathsf{Peak}(\eta_i) = \mathcal{O}_i \left(1 + \sum_{k\in\mathcal{A}}\nolimits \mathcal{O}_k\right)^{\frac{\sum_{k\in\mathcal{A}}\nolimits \mathcal{O}_k}{1+\sum_{k\in\mathcal{A}}\nolimits \mathcal{O}_k}}.
\]
To interpret this formula, consider the specific situation where all areas \emph{other} than area $i$ have perfect bias tunings. Then $\mathcal{O}_k = 0$ for all $k \neq i$, and the peak response simplifies to 
\[
\mathsf{Peak}(\eta_i) = \mathcal{O}_i \left(1 + \mathcal{O}_i\right)^{\frac{\mathcal{O}_i}{1+\mathcal{O}_i}}.
\]
This function is plotted in Figure \ref{Fig:OverUnder} for the range $\mathcal{O}_i \in [-\tfrac{1}{2},\tfrac{1}{2}]$.

Observe that $\mathsf{Peak}(\eta_i) \approx \mathcal{O}_i$ in a large range around $\mathcal{O}_i = 0$. It follows that for realistic errors in the bias tuning, underbiasing does not cause power withdrawl appreciably greater in magnitude than the corresponding power support provided by overbiasing. For very large bias tuning errors however, withdrawl is indeed greater in magnitude than the corresponding support. Therefore, Cohn's second claim is true, but only for very extreme overbiased and underbiased tunings.

\smallskip

\begin{figure}[ht!]
\begin{center}
\includegraphics[width=0.95\columnwidth]{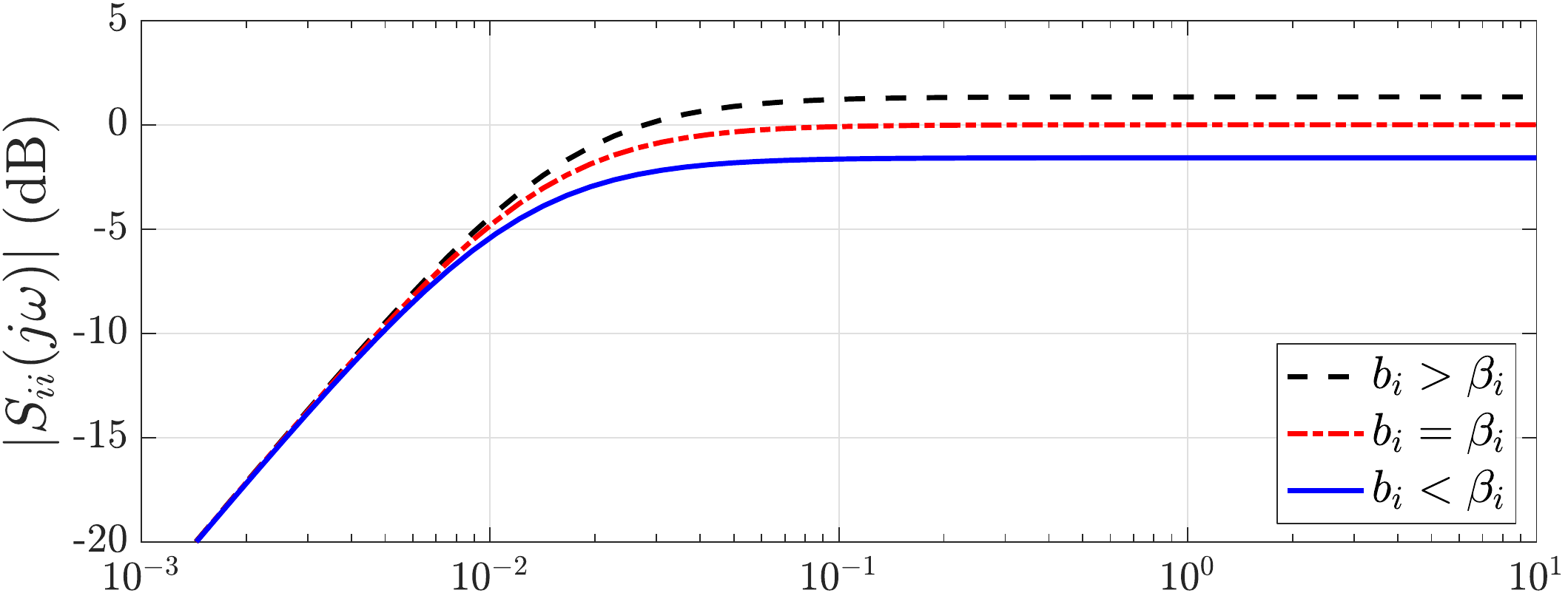}
\caption{Bode plot of $S_{ii}(s)$ for overbiased and underbiased tunings.}
\label{Fig:Bode}
\end{center}
\end{figure}

\emph{Transfer function $S_{ij}(s)$ from net load disturbance $\Delta P_{j}^{\rm L}$ in area $j$ to area control error $\mathsf{ACE}_i$ in area $i$:} Very similar calculations using \eqref{Eq:SlowTimeScaleAGC2-LTI} show that
\[
\begin{aligned}
S_{ij}(s) 
&= -\frac{\tau^{\prime}s}{\tau^{\prime}s+1}\left[\delta_{ij} - \frac{1}{\beta}(\beta_i-b_i)\frac{\tau^{\prime}s}{\tau^{\prime}s+\tfrac{1}{\beta}\sum_{k}b_k}\right].
\end{aligned}
\]
A representative Bode plot of $S_{ii}(s)$ for underbiased and overbiased tunings is shown in Figure \ref{Fig:Bode}. The peak sensitivity can be computed to be
\begin{equation}\label{Eq:HInfNorm}
\sup_{\omega \in \real} |S_{ii}(j\omega)| = |1 - \tfrac{1}{\beta}(\beta_i-b_i)| = |1 + \mathcal{O}_i|.
\end{equation}

%

Overbiasing results in aggressive AGC response to frequnecy deviations, and tends to \emph{increase} the high-frequency sensitivity of the local ACE with respect to local disturbances; moderate underbiasing has the opposite effect. We conclude that there is a natural dynamic tension in bias tuning, between providing proper assistance to other areas (Figure \ref{Fig:OverUnderStep}) and reducing local sensitivity of the ACE to load changes (Figure \ref{Fig:Bode}). As a final word of caution, it is important to note that when the bias is tuned imperfectly, the numerical value of the ACE \emph{does not equal} the generation-load mismatch; it is merely a proxy (Lemma \ref{Lem:ACE}). See \cite{JWSP:20e} for further comments on this gap.

\section{Simulations on Two-Area Test System}
\label{Sec:Simulations}

\begin{figure}[ht!]
\begin{center}
\includegraphics[width=0.95\columnwidth]{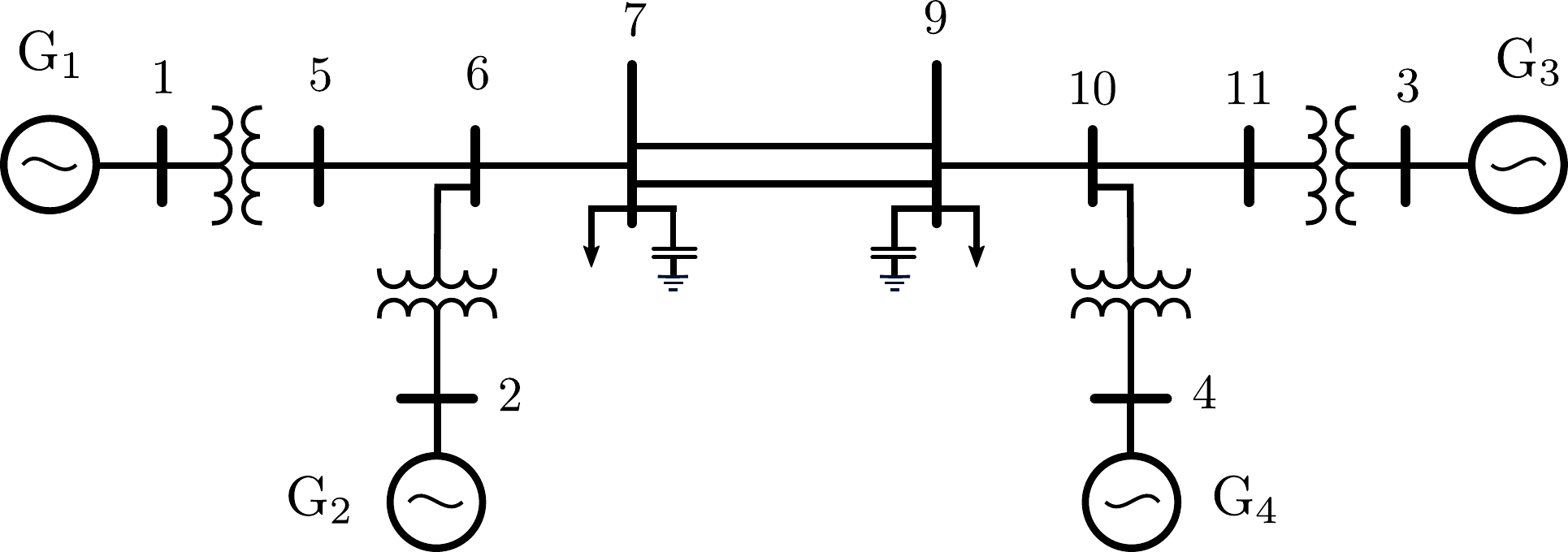}
\caption{Two-area four-machine test system.}
\label{Fig:KundurLine}
\end{center}
\end{figure}

We illustrate our results with straightforward simulations of the AGC controller \eqref{Eq:AGCSimple} on the Kundur two-area four-machine test system (Figure \ref{Fig:KundurLine}), implemented in MATLAB's Simscape Electrical. The system is three-phase and includes full-order machine, turbine-governor, excitation, and PSS models; SVCs were integrated at buses 7 and 9 to support voltage levels.

\begin{figure}[ht!]
\centering
\begin{subfigure}{0.99\linewidth}
\includegraphics[width=\linewidth]{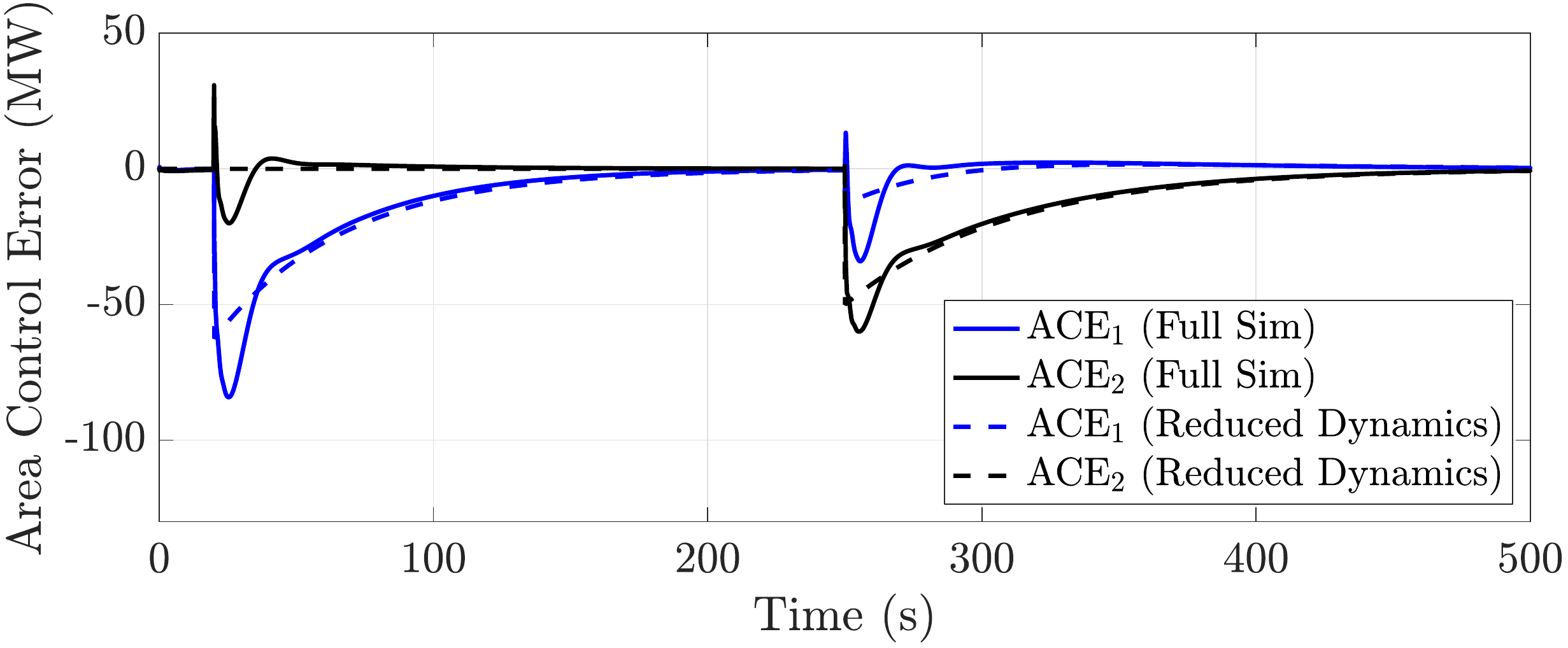}
\caption{Area control errors.}
\label{Fig:Kundur-a}
\end{subfigure}\\
\smallskip
\begin{subfigure}{0.99\linewidth}
\includegraphics[width=\linewidth]{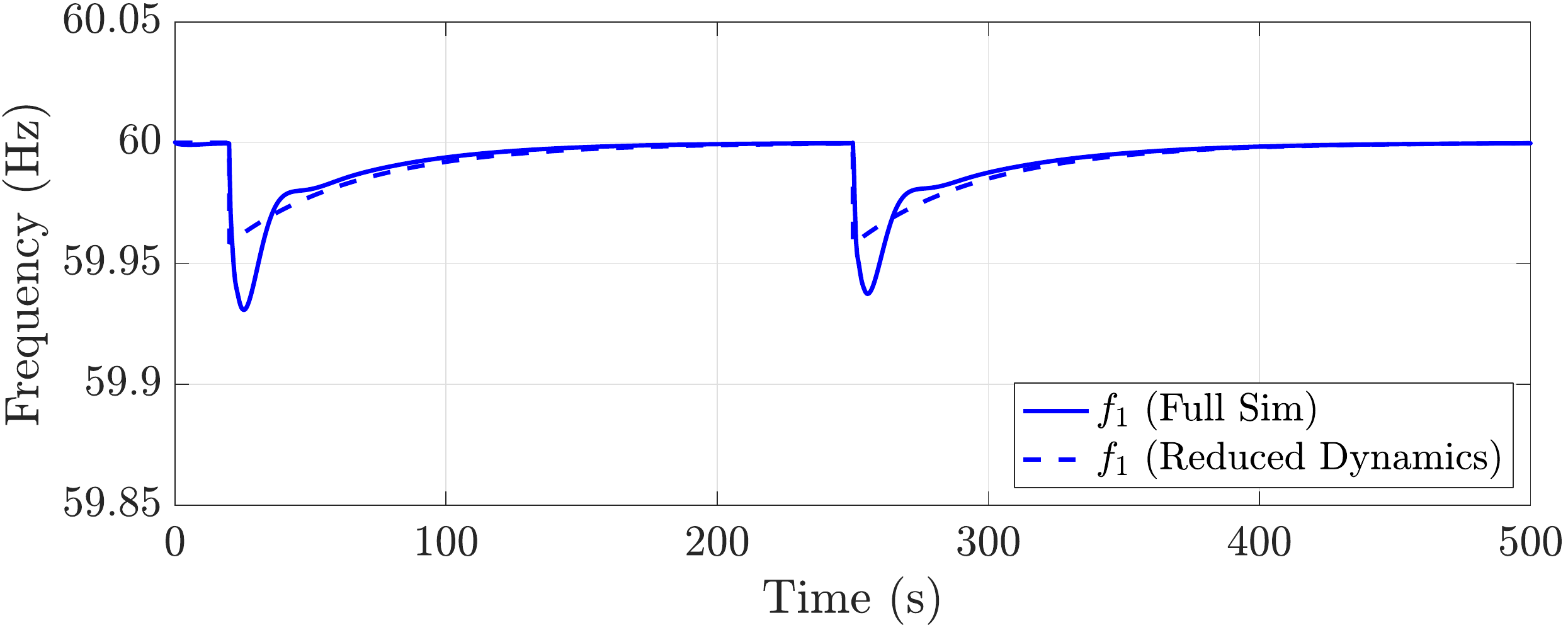}
\caption{Frequency in area 1 ($\approx$ frequency in area 2).}
\label{Fig:Kundur-b}
\end{subfigure}\\
\smallskip
%
%
%
\smallskip
\begin{subfigure}{0.99\linewidth}
\includegraphics[width=\linewidth]{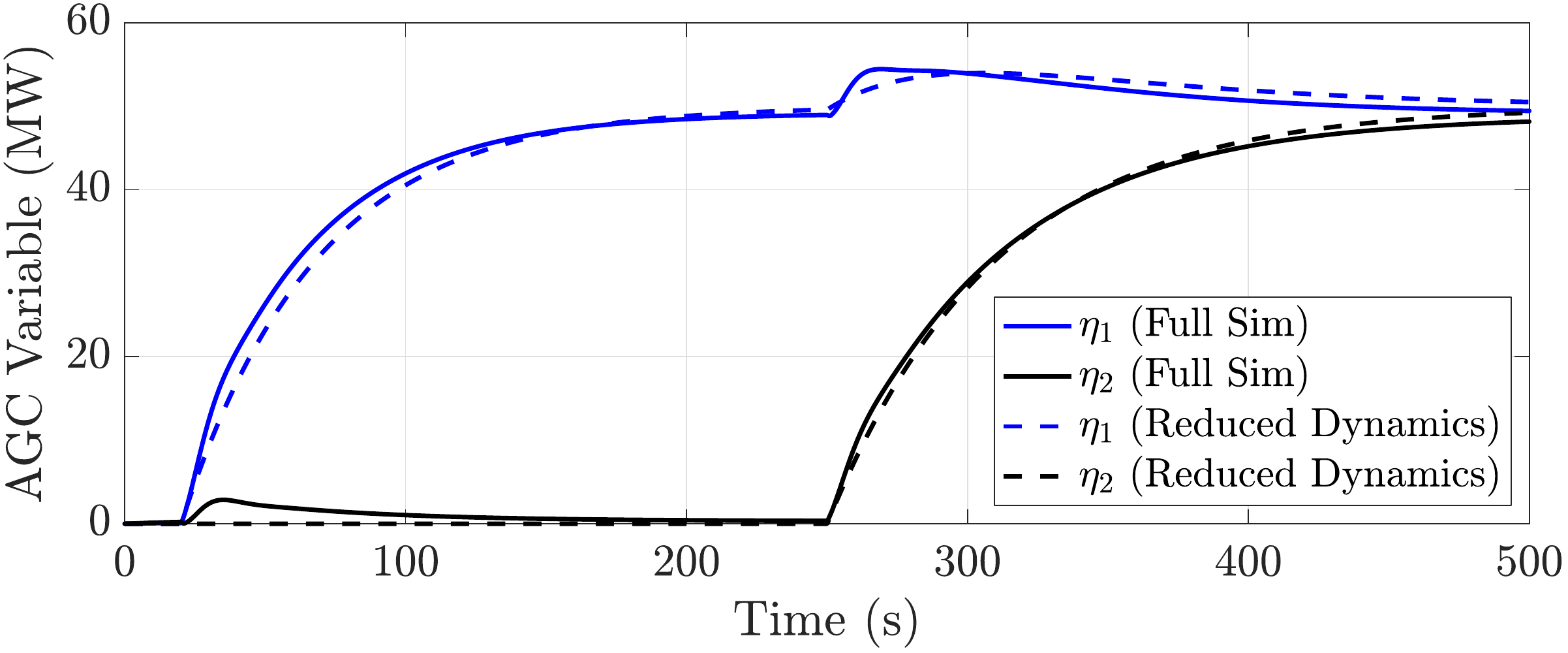}
\caption{AGC control variables.}
\label{Fig:Kundur-e}
\end{subfigure}
\caption{Response of two-area four-machine system with overbiasing in area 1 controller.}
\label{Fig:Kundur}
\end{figure}

 All four generators have 5\% primary governor droop, and there is no appreciable load-frequency damping in the system; it follows that $\beta_1 = \beta_2 = 40$ p.u./p.u. The AGC controller in area 2 is tuned correctly with $b_2 = \beta_2$, while the controller in area 1 is overbiased with $b_1 = 1.5\beta_1$. The time constants are $\tau_1 = \tau_2 = 60$s. Only generators G1 and G3 participate in AGC. The system is subject to a 50MW load increase at bus 7 in area 1 at $t = 20$s, and then a 50MW load increase at bus 9 in area 2 at $t = 250$s. 
We compare the full-order simulation results with those obtained from the reduced dynamics \eqref{Eq:SlowTimeScaleAGC2} by plotting in Figure \ref{Fig:Kundur} the resulting ACEs, the system frequency, and the control variables $\eta_i$. All plots demonstrate that \textemdash{} aside from transients after the disturbances due to primary control \textemdash{} the reduced dynamics \eqref{Eq:SlowTimeScaleAGC2} provide an excellent approximation of the full nonlinear response. 


\begin{figure}[ht!]
\centering
\includegraphics[width=\linewidth]{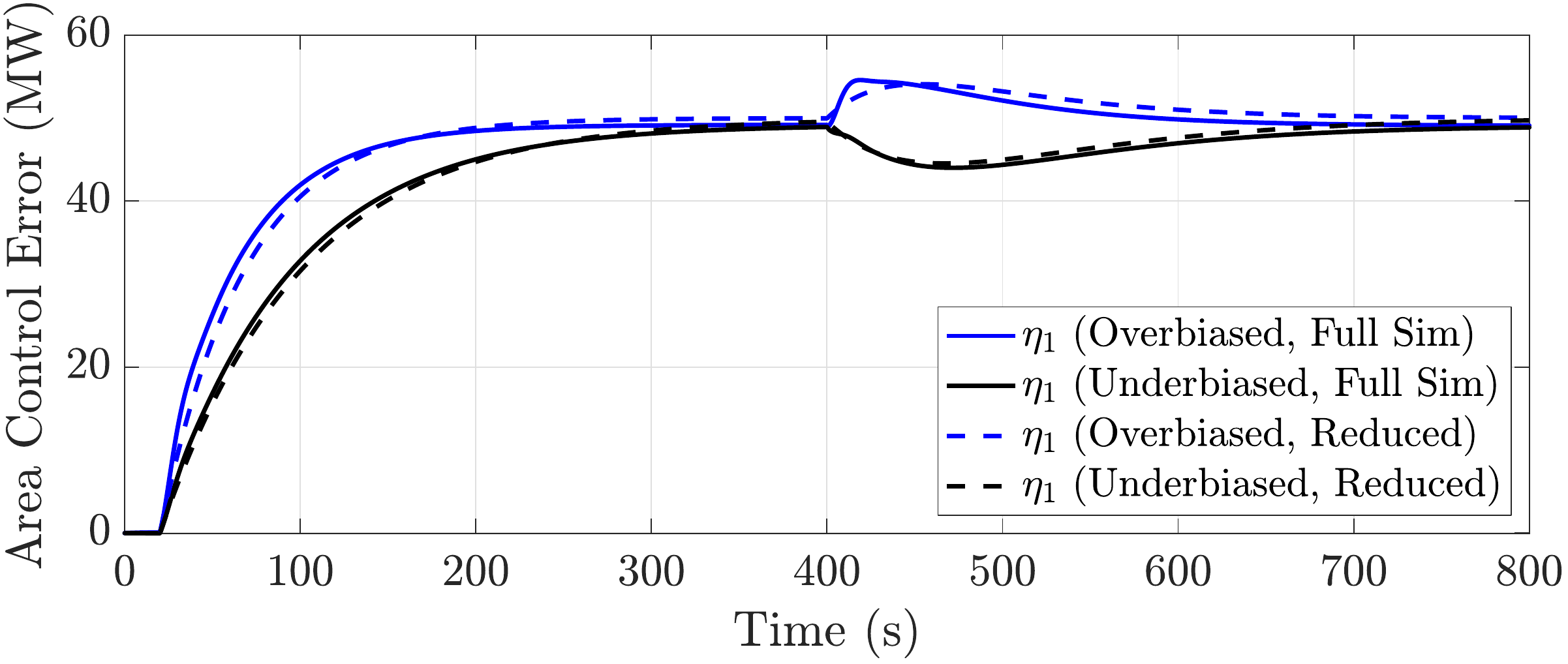}
\caption{Comparison of $\eta_1$ for overbiased and underbiased tunings.}
\label{Fig:KundurOverUnder}
\end{figure}

In Figure \ref{Fig:KundurOverUnder} we repeat the same disturbance, and compare the previous results with what one would obtain with an \emph{underbiased} tuning $b_1 = 0.5\beta_1$ in area 1; in our notation of Section \ref{Sec:DynamicPerf}, we have $\mathcal{O}_1 = \pm 0.25$ for the two tunings. First, notice that the underbiased response is slower, which is consistent with our eigenvalue result from Section \ref{Sec:DynamicPerf}. Second, note that Figure \ref{Fig:KundurOverUnder} is consistent with Figure \ref{Fig:OverUnder}, as the peak responses in area 1 to the disturbance in area 2 are roughly equal and opposite for overbiased vs. underbiased tunings. 

\section{Discussion}
\label{Sec:Discussion}

A possible objection to our methodology is that the turbine/governor and energy supply (e.g., boiler) dynamics of some traditional power plants are \emph{so} slow that their interaction with AGC dynamics should not be neglected (e.g., \cite{CWT-RLC:76,GZ-JM-QW:19}). It would appear difficult to make any general statements concerning this which are broadly applicable across all power systems. If boiler dynamics are to be retained, then the preceding analysis must be modified; this is a straightforward extension for future work. Another possible objection is that dynamically important nonlinear physical effects such as turbine-governor deadbands, and digital implementation effects such as discrete sampling, have not been included in the analysis. These effects further degrade system performance by limiting controller sensitivity and bandwidth. Our results here identify the stability and performance characteristics of ``best case'' AGC implementations, and therefore reveal the performance characteristics of AGC which are intrinsic to the decentralized control structure, the measurements used, and the basic physics of power systems. Our work here is thus a pre-requisite for rigorously understanding any additional performance-degrading dynamic effects.


%

\section{Conclusions}
\label{Sec:Conclusions}

We have presented a rigorous nonlinear dynamic analysis of AGC in interconnected power systems. Our approach provides a simple methodology for understanding the system-wide dynamics induced by AGC on a long time-scale, and we have exploited this to understand some of the fundamental stability and performance properties of AGC. Simulation results verify that the reduced dynamics \eqref{Eq:SlowTimeScaleAGC2} quite accurately model the slow time-scale dynamic behaviour of AGC, and can therefore be used to quantify dynamic performance. Among other conclusions, our results provide rigorous control-theoretic justification for overbiasing of AGC systems in practice.

There are many open avenues for extensions of this analysis, some of which have already been noted, and one of which is explored in \cite{JWSP:20e}. Another interesting direction concerns so-called doubly-integrated AGC schemes, and more broadly, the implications of the current analysis for management of inadvertent interchange between control areas.

\appendices
\section{Technical Results and Proofs}
\label{App:1}

\begin{pfof}{Lemma \ref{Lem:ACE}}
\ref{Lem:ACE-2} $\Longrightarrow$ \ref{Lem:ACE-3}: This is immediate.

\ref{Lem:ACE-3} $\Longrightarrow$ \ref{Lem:ACE-2}: In synchronous steady-state, we have from \eqref{Eq:Synchronous} that $\Delta f_1 = \cdots = \Delta f_N = \Delta f_{\rm ss}$, and therefore
\begin{equation}\label{Eq:ACEZeroing}
\begin{aligned}
0 = \mathsf{ACE}_k &= \Delta \NI_{k} + b_k\Delta f_{\rm ss}
\end{aligned}
\end{equation}
for all $k \in \mathcal{A}$. Summing \eqref{Eq:ACEZeroing} over all areas $k$ and using \eqref{Eq:NetBalance}, we find that $0 = (b_1+\cdots+b_N)\Delta f_{\rm ss}$ which implies that $\Delta f_{\rm ss} = 0$; it now follows from \eqref{Eq:ACEZeroing} that $\Delta \NI_{k} = 0$ for all $k$.

\ref{Lem:ACE-1} $\Longleftrightarrow$ \ref{Lem:ACE-3}: Substituting \eqref{Eq:GovPower} into \eqref{Eq:Balancek} and rearranging, we obtain
\[
\sum_{i\in\mathcal{G}_k}\nolimits (u_{k,i}-u_{k,i}^{\star}) = \Delta P_{k}^{\rm L} + \Delta p_{k}^{\rm tie} +  \beta_k \Delta f_{k}
\]
Using \eqref{Eq:PowerFlowNI} on the right-hand side and \eqref{Eq:GovIdentity} on the left-hand side, we find that
\[
\sum_{i\in\mathcal{G}_k^{\rm AGC}}\nolimits (u_{k,i}-u_{k,i}^{\star}) = \Delta P_{k}^{\rm L} + \Delta \NI_{k} +  \beta_k \Delta f_{k}.
\]
Using the same argument as in \ref{Lem:ACE-3} $\Longleftrightarrow$ \ref{Lem:ACE-2}, it is straightforward to argue that $\Delta \NI_{k} +  \beta_k \Delta f_{k} = 0$ for all $k \in \mathcal{A}$ if and only if $\mathsf{ACE}_k = 0$ for all $k \in \mathcal{A}$, and the equivalence \ref{Lem:ACE-1} $\Longleftrightarrow$ \ref{Lem:ACE-3} therefore follows.
\end{pfof}

\begin{lemma}[\bf Interconnection Matrix]\label{Lem:Inter}
Let $\boldsymbol{b},\boldsymbol{\beta} \in \real^N$ have strictly positive elements with $\beta = \vones[N]^{\T}\boldsymbol{\beta}$, and set $\boldsymbol{\gamma} = (\boldsymbol{\beta}-\boldsymbol{b})/\beta$. Consider the matrix $\mathcal{B} = I_N - \boldsymbol{\gamma}\vones[N]^{\T}$. Then 
\begin{enumerate}
\item \label{Lem:Inter-1} $\mathrm{eig}(\mathcal{B}) = \{1,1,\ldots,1,(1-\vones[]^{\T}\boldsymbol{\gamma})\}$,
\item \label{Lem:Inter-2} $-\mathcal{B}$ is diagonally stable,
\item \label{Lem:Inter-3} $\mathcal{B}^{-1} = I_N + \tfrac{1}{1-\vones[N]^{\T}\boldsymbol{\gamma}}\boldsymbol{\gamma}\vones[N]^{\T}$, and
\item \label{Lem:Inter-4} $(sI_{N}+\mathcal{B})^{-1} = \frac{1}{s+1}\left[I_N + \tfrac{1}{(s+1) -\vones[]^{\T}\boldsymbol{\gamma}}\boldsymbol{\gamma}\vones[N]^{\T}\right]$.
\end{enumerate}
\end{lemma}

\begin{pfof}{Lemma \ref{Lem:Inter}}
Items \ref{Lem:Inter-1},\ref{Lem:Inter-3},\ref{Lem:Inter-4} are by direct (if somewhat tedious) calculation. Item \ref{Lem:Inter-2} follows directly from \cite[Theorem 2.1]{NM-JWSP:20m}.
\end{pfof}

\begin{pfof}{Theorem \ref{Thm:AGCStable}}
For the closed-loop system \eqref{Eq:NonlinearPowerSystem} and \eqref{Eq:AGC},\eqref{Eq:Allocation} let $\varepsilon > 0$ be small and define $\tilde{\tau}_k \define \varepsilon\tau_k$ for some values $\tilde{\tau_k} > 0$ which are $\mathcal{O}(1)$. Defining the new time variable $\ell = \varepsilon t$ leads to the singularly perturbed system \cite[Chp. 11]{HKK:02}
\begin{subequations}\label{Eq:SingularPerturbation}
\begin{align}
\label{Eq:SingularPerturbation-1}
\varepsilon \frac{\mathrm{d}x}{\mathrm{d}\ell} &= F(x,u,w)\\
\label{Eq:SingularPerturbation-2}
\tilde{\tau}_k\frac{\mathrm{d}\eta_k}{\mathrm{d}\ell} &= -(\Delta \NI_{k} + b_k\Delta f_{k}) - \vones[]^{\T}(u_k-P_k)\\
\label{Eq:SingularPerturbation-3}
u_{k,i} &= \mathrm{sat}_{k,i}(u_{k,i}^{\star} + \alpha_{k,i}\eta_k),
\end{align}
\end{subequations}
for $k \in \mathcal{A}$ and $i \in \mathcal{G}_{k}^{\rm AGC}$, to which we will apply Lyapunov arguments (e.g., \cite{AS-HKK:84}, \cite[Theorem 11.3]{HKK:02}). The ``boundary layer'' dynamics are \eqref{Eq:SingularPerturbation-1} with $u$ considered as a parameter; Assumptions \ref{Ass:PowerSystem-1}--\ref{Ass:PowerSystem-3} guarantee that the conditions imposed in \cite[Lemma 1]{AS-HKK:84} are satisfied. Using Assumption \ref{Ass:PowerSystem-4} and \ref{Ass:PowerSystem-5}, the calculations preceding \eqref{Eq:SlowTimeScaleAGC} have shown that, modulo time-scaling, the reduced dynamics associated with \eqref{Eq:SingularPerturbation} are exactly given by the nonlinear system \eqref{Eq:SlowTimeScaleAGC}. We now argue that the reduced dynamics possess an equilibrium point which is both globally asymptotically stable and locally exponentially stable. By Assumption \ref{Ass:Feasibility} we have that $\Delta P^{\rm L}_{k} \in \mathcal{C}_k$. Examining the definition of $\varphi_k(\eta_k)$, one can see that $\varphi_{k}(\eta_k)$ is a piecewise linear and non-decreasing function of $\eta_k \in \real$ and that $\mathrm{image}(\varphi_k) = \mathrm{closure}(\mathcal{C}_k)$. The preimage $\mathcal{P}_k = \setdef{\eta_k \in \real}{\varphi_{k}(\eta_k) \in \mathcal{C}_k}$ of $\mathcal{C}_k$ can be  explicitly computed to be the interval
\[
\mathcal{P}_k \define \left(\min_{i\in\mathcal{G}_k}\tfrac{1}{\alpha_i}(\underline{u}_{k,i}-u_{k,i}^{\star}),\max_{i\in\mathcal{G}_k}\tfrac{1}{\alpha_i}(\overline{u}_{k,i}-u_{k,i}^{\star})\right),
\]
and a simple argument shows that $\varphi_k$ is a \emph{strictly} increasing function on $\mathcal{P}_k$. We conclude that $\map{\varphi_{k}}{\mathcal{P}_k}{\mathcal{C}_k}$ is a bijective function, and hence there exists a unique $\bar{\eta}_k \in \mathcal{P}_k$ such that $\varphi_{k}(\bar{\eta}_k) = \Delta P^{\rm L}_{k}$ for all $k \in \mathcal{A}$. Since $\mathcal{B}_{\rm txt}$ is invertible (Lemma \ref{Lem:Inter}), we conclude that $\bar{\eta} = (\bar{\eta}_1,\ldots,\bar{\eta}_N)$ is the unique equilibrium point of \eqref{Eq:SlowTimeScaleAGC}. By Lemma \ref{Lem:Inter}, the matrix $\mathcal{B}_{\rm txt}$ is diagonally stable, so there exists a diagonal matrix $D = \mathrm{diag}(d_{1},\ldots,d_N) \succ 0$ such that $Q \define \mathcal{B}_{\rm txt}^{\T}D + D\mathcal{B}_{\rm txt} \succ 0$. For \eqref{Eq:SlowTimeScaleAGC}, define the Lyapunov candidate
\[
V(\eta) = \tfrac{1}{2}\sum_{k=1}^{N}\nolimits d_k \tau_k\int_{\bar{\eta}_k}^{\eta_{k}} (\varphi_k(\xi_k) - \varphi_k(\bar{\eta}_k))\,\mathrm{d}\xi_k,
\]
which obviously satisfies $V(\bar{\eta}) = 0$. Since $\varphi_k$ is non-increasing and is strictly increasing on $\mathcal{P}_{k}$, we have that $V(\eta) > 0$ for all $\eta \neq \bar{\eta}$, and that $V$ is radially unbounded. We compute along trajectories of \eqref{Eq:SlowTimeScaleAGC} that
\[
\begin{aligned}
\dot{V}(\eta) &= \tfrac{1}{2}(\varphi(\eta)-\varphi(\bar{\eta}))^{\T}D(-\mathcal{B}_{\rm txt})(\varphi(\eta)-\Delta P^{\rm L}_{k})\\
&= \tfrac{1}{2}(\varphi(\eta)-\varphi(\bar{\eta}))^{\T}D(-\mathcal{B}_{\rm txt})(\varphi(\eta)-\varphi(\bar{\eta}))\\
&= -(\varphi(\eta)-\varphi(\bar{\eta}))^{\T}Q(\varphi(\eta)-\varphi(\bar{\eta})) < 0
\end{aligned}
\]
for all $\eta \neq \bar{\eta}$. We conclude that $\bar{\eta}$ is globally asymptotically stable. Since $\bar{\eta}_k \in \mathcal{P}_k$ for all $k \in \mathcal{A}$ and $\mathcal{P}_{k}$ is an open interval on which $\varphi_k$ is piecewise linear and strictly increasing, there exists some $r_k > 0$ such that $\varphi_k$ is both strongly monotone and Lipschitz continuous for all $\eta_k \in [\bar{\eta}_k - r_k,\bar{\eta}_k + r_k]$. Combining this with the above calculations, it follows quickly that $V$ is \emph{locally} a so-called \emph{quadratic-type Lyapunov function}, and that $\bar{\eta}$ is therefore locally exponentially stable. The remaining conditions of \cite[Lemma 1]{AS-HKK:84} are now satisfied, and it now follows that the equilibrium $(\bar{x},\bar{\eta})$ of \eqref{Eq:SingularPerturbation} \textemdash{} and hence, of the closed-loop system \textemdash{} is locally exponentially stable for sufficiently small $\varepsilon > 0$, or equivalently, for sufficiently large $\min_{k} \tau_k$. Finally, since $\mathcal{B}_{\rm ACE}$ is invertible (Lemma \ref{Lem:Inter}), we see from \eqref{Eq:SlowTimeScaleAGC} that $\mathsf{ACE}_k = 0$ at equilibrium, which completes the proof. The proof for the simplified AGC controller \eqref{Eq:AGCSimple} is nearly identical.
\end{pfof}

\IEEEpeerreviewmaketitle



\ifCLASSOPTIONcaptionsoff
  \newpage
\fi


\renewcommand{\baselinestretch}{1}
\bibliographystyle{IEEEtran}

\bibliography{/Users/jwsimpso/GoogleDrive/JohnSVN/bib/brevalias,%
/Users/jwsimpso/GoogleDrive/JohnSVN/bib/Main,%
/Users/jwsimpso/GoogleDrive/JohnSVN/bib/JWSP,%
/Users/jwsimpso/GoogleDrive/JohnSVN/bib/New%
}

\begin{IEEEbiography}[{\includegraphics[width=1in,height=1.25in,clip,keepaspectratio]{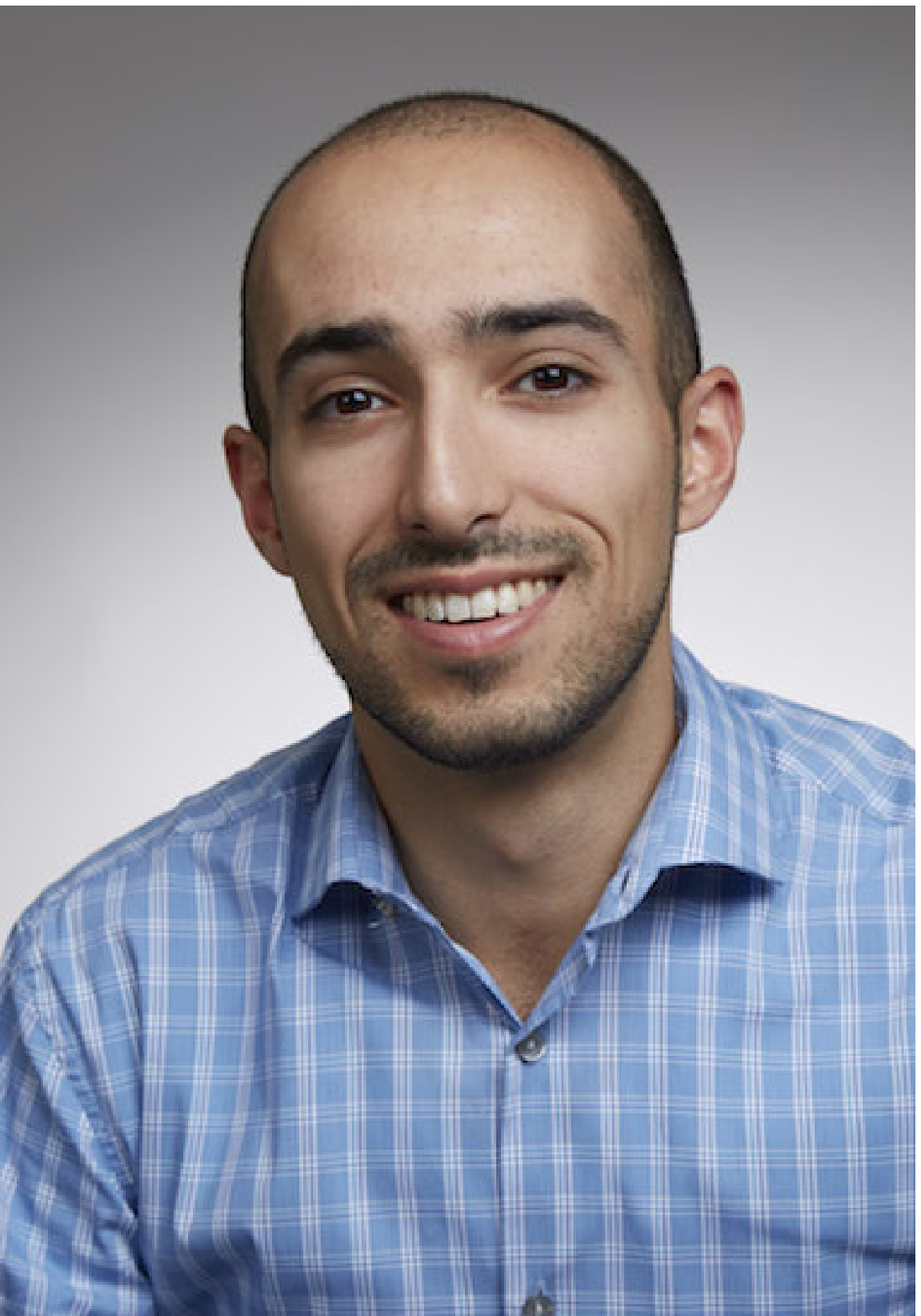}}]{John W. Simpson-Porco} (S'11--M'16) received the B.Sc. degree in engineering physics from Queen's University, Kingston, ON, Canada in 2010, and the Ph.D. degree in mechanical engineering from the University of California at Santa Barbara, Santa Barbara, CA, USA in 2015.

He is currently an Assistant Professor of Electrical and Computer Engineering at the University of Toronto, Toronto, ON, Canada. He was previously an Assistant Professor at the University of Waterloo, Waterloo, ON, Canada and a visiting scientist with the Automatic Control Laboratory at ETH Z\"{u}rich, Z\"{u}rich, Switzerland. His research focuses on feedback control theory and applications of control in modernized power grids.

Prof. Simpson-Porco is a recipient of the 2012--2014 IFAC Automatica Prize and the Center for Control, Dynamical Systems and Computation Best Thesis Award and Outstanding Scholar Fellowship. He currently serves as an Associate Editor for IEEE Transactions on Smart Grid.
\end{IEEEbiography}

\end{document}